\documentclass{amsart}
\usepackage{amsmath}
\newtheorem{theorem}{Theorem}[section]

\theoremstyle{definition}

\theoremstyle{remark}
\newtheorem{remark}[theorem]{Remark}

\numberwithin{equation}{section}

\newcommand{\abs}[1]{\lvert#1\rvert}
\newcommand{\bs}{\bigskip}

\newcommand{\s}{\sigma}
\newcommand{\T}{\tau}

\begin{document}

\title[An asymptotically Hilbertian space which fails the A.P.]{An
example of an asymptotically Hilbertian space  which fails the
approximation property}

\author{P.~G. ~Casazza}
\address{Department of Mathematics, University of
Missouri-Columbia, Columbia, MO 65211 U.S.A.}
\email{pete@math.missouri.edu}
\thanks{The first author was supported by NSF grant DMS-970618.}

\author{C.~L. ~Garc\'{\i}a}
\address{Department of Mathematics, Texas A\&M University
College Station, TX  77843--3368 U.S.A}
\email{clgarcia@math.tamu.edu}

\author{W.~B. ~Johnson}
\address{Department of Mathematics, Texas A\&M University
College Station, TX  77843--3368 U.S.A}
\email{johnson@math.tamu.edu}
\thanks{The second and third authors were supported in part by
NSF grant DMS-9623260, DMS-9900185, and by Texas Advanced
Research Program under Grant No. 010366-163.}


\date{March 1, 2000.}

\commby{Nicole Tomczak-Jaegermann}


\begin{abstract}
Following Davie's example of a Banach space
failing the approximation property (\cite{D}), we show  how to
construct a Banach  space $E$ which is asymptotically Hilbertian
and fails the approximation  property. Moreover, the space $E$ is
shown to be a subspace
of a space with an unconditional basis which is  ``almost''  a weak
Hilbert  space and which can be written as the direct sum of two
subspaces
all of whose subspaces have the approximation property.
\end{abstract}

\maketitle

\section{Introduction}

This paper is concerned with the relationship between the approximation 
property and notions about Banach spaces which are in some sense close to Hilbert space, namely, the notion of asymptotically Hilbertian space and of 
weak Hilbert space.

The spaces we discuss are of the form 
$Z = \left(\sum_{n=0}^{\infty}\ell_{p_n}^{k_n}\right)_{2}$
with $p_n\downarrow 2$ and $k_n\uparrow \infty$. 
It is easy to check that any space of this form is an 
asymptotically Hilbertian space (see below for definitions). For particular
sequences $(p_n)$ and $(k_n)$ we show that such a $Z$ has a subspace $E$ failing the approximation property. Moreover, 
we can choose a subsequence of $(p_n)$, such that if
$N_1= \{j  | p_{n_{2k+1}} \leq j < p_{n_{2(k+1)}}, k\geq 0\}$  
and $N_2 = \mathbb{N} - N_1$ then for
$Z_i = \left(\sum_{j\in N_i}\ell_{p_j}^{k_j}\right)_{2}, i=1,2$, we have that
$Z = Z_1\oplus Z_2$ and that all subspaces of $Z_1$ and of $Z_2$ have the approximation property (\cite{J}).

The construction of $E$ provides quantitative estimates which
show that $Z$ and hence also $E$ is surprisingly close to being
a weak Hilbert space (note that weak Hilbert spaces enjoy the approximation property \cite{P}).  

First we recall the notion of \textsl{asymptotically Hilbertian space}. Given  integers $n\ge 0$, $m\ge 1$  and a constant $K$, say that
$X$ satisfies
$H(n,m,K)$ provided there is an $n$-codimensional subspace $X_m$ of $X$
so 
that every $m$-dimensional subspace of $X_m$ is $K$-isomorphic to
$\ell_2^m$.  A Banach space $X$ is
said to be \textsl{asymptotically Hilbertian} provided there is a
constant
$K$ so that for every $m$ there exists $n$ so that $X$ satisfies
$H(n,m,K)$.  Since here we are interested in good estimates, we denote
by
$H_X(m,K)$ the smallest $n$ for which $X$ has $H(n,m,K)$.  Thus if $X$
is
$K$-isomorphic to a Hilbert space, then $H_X(m,K)=0$ for all $m$.  The
growth rate of $H_X(m,K)$ for a fixed $K$ as $m\to\infty$ is one
measurement of the closeness of $X$ to a Hilbert space.  

A Banach space is called a \textsl{weak Hilbert space\/} provided that
there are positive constants $\delta$ and $K$ so that for every $n$,
every $n$ dimensional subspace of $X$ contains a further subspace $E$ of
dimension at least $\delta n$ so that $E$ is $K$-isomorphic to a Hilbert
space and $E$ is $K$-complemented in $X$ (that is, there is a projection
having norm at most $K$ from $X$ onto $E$).  

The definition  of the property $H(n,m,K)$ was made in \cite{J},
although the nomenclature ``asymptotically Hilbertian" was coined by
Pisier \cite{P}. Weak Hilbert spaces were introduced by Pisier
\cite{P}, who gave many equivalences to the property of being weak
Hilbert; we chose the one most relevant for this paper as the
definition of weak Hilbert.  

Relations between the weak Hilbert property and the asymptotically  
Hilbertian property are given in \cite{J} and \cite{P}.   First, a weak
Hilbert space must be asymptotically Hilbertian (\cite[Section 4]{P}).
It seems likely that if $X$ is weak Hilbert then for some
$K$,
$H_X(m,K)\le K m$, but in fact no reasonable estimates are known for
$H_X(m,K)$ when
$X$ is a weak Hilbert space.   It is known (see \cite{CJT}, \cite{NT-J})
that if
$X$ is a weak Hilbert space which has an unconditional basis and there
is a
$K$ so that $H_X(m,K)$ is dominated by $f(m) $ for some iterate $f$ of
$\exp$, then for any iterate $g$ of $\log$ there is another constant
$K'$
so that 
$H_X(m,K')\le K'g(m)$. In the other direction, it follows from \cite{J}
that if
for some $K$ the sequence $H_X(m,K)$ grows sufficiently slowly as
$m\to \infty$ (say, like   $\log\log m$), then $X$ is a weak
Hilbert space.  In this paper we are interested in examples of spaces
which are of type $2$. We refer to Chapter 11 of \cite{DJT} or Section
1.4 of  \cite{T-J} for the definitions and basic theory of type $p$ and
cotype
$p$ as well and the type $p$ and cotype $p$ constants $T_p(X)$,
$C_p(X)$ of a Banach space $X$.  Relevant for us is that if
$X$ is a type
$2$ space and
$E$ is a subspace of $X$ which is
$K$-isomorphic to a Hilbert space then, by Maurey's extension theorem,
$E$ is
$T_2(X)K$- complemented in $X$ (\cite[Corollary 12.24]{DJT}). Thus
it is clear that if
$X$ is of type
$2$ and for some $K$, $H_X(m,K)\le K m$, then $X$ is weak Hilbert.  Here
we should mention that by \cite{FLM}, polynomial   growth of
$H_X(m,K)$ (as $m\to \infty$)  implies linear growth of $H_X(m,K')$ for
some
$K'$.

Our main interest here is the linkage among the weak Hilbert
property, the asymptotically Hilbertian property, and the approximation
property.   The arguments in \cite{J} show that if $X$ has type $2$ and 
for some $K$, $H_X(m,K)\le K\log m$ for infinitely many $m$, then all
subspaces of $X$ (even all subspaces of every quotient of $X$) have the
approximation property.  In \cite{P} it is shown that all weak Hilbert
spaces have the approximation property.  Thus if
$X$ is of type $2$ and for some $K$, $H_X(m,K)\le K m$ for {\sl all\/}
$m$, then   all subspaces of every quotient of $X$ have the
approximation
property. It is easy to build examples of a  type $2$ space $X$ for
which
there is a constant $K$ so that for any iterate $f$ of the $\log$
function, 
$H_X(m,K)\le Kf(m)$ for infinitely many $m$ and yet $X$ is not a weak
Hilbert space.  Now such a space $X$ is in some sense close  to Hilbert
space and, in particular, every subspace of $X$ has the approximation
property.  In this paper we show that there are two such spaces;  call
them $Z_1$ and $Z_2$; so that $Z:=Z_1\oplus Z_2$ has a subspace which
fails the approximation property. Moreover, $Z$ has an unconditional
basis ($Z=\sum_{n=0}^\infty \ell_{p_n}^{k_n}$ for appropriate
$p_n\downarrow 2$ and $k_n\uparrow \infty$) and is nearly a weak Hilbert
space in the sense that for some
$K$, the growth rate of
$H_Z(m,K)$ as
$m\to \infty$ is close to being polynomial in $m$ ($H_Z(m,K)\le m^{\log
\log m}$ is what we get; recall that polynomial  growth 
of
$H_X(m,K)$ gives linear growth of $H_X(m,K')$ for some $K'$.).

\section{the example}

We follow closely A.~M. ~Davie's construction of a Banach space
failing the approximation property \cite{D}. Davie constructed for $p > 2$
a subspace of $\ell_p = \left(\sum \ell_{p}^{k_n}\right)_{p}$ which fails
the approximation property. He could as well have used
$\left(\sum \ell_{p}^{k_n}\right)_{r}$ for any $1 \leq r \leq \infty$. Here
we use instead $Z := \left(\sum_{n=0}^{\infty}\ell_{p_n}^{k_n}\right)_{2}$
where $p_n\downarrow 2$ appropriately and $k_n$ as in \cite{D}. Basically we
compute how fast $p_n$ can go to $2$ so that Davie's argument yields a 
subspace of $Z$ which fails the a.p.. The obvious condition is that
$k_n^{1/2 -1/{p_n}}$ cannot be bounded, for if $k_n^{1/2 -1/{p_n}}$ is
bounded then $\left(\sum_{n=0}^{\infty}\ell_{p_n}^{k_n}\right)_{2}$ is isomorphic to $\ell_2$.
  
For any integer $n \ge 0$ consider an Abelian group $G_n$ of order
$k_n = 3\cdot 2^n$, and let $\s_1^n, \ldots, \s_{2^n}^n, \T_1^n,
\ldots, \T_{2^{n+1}}^n$ be the characters of $G_n$. Lemma $(b)$ in
\cite{D} shows that this enumeration of the characters of $G_n$ can be
chosen so  that  there exists an absolute constant $A > 0$ such that
for all $g \in G_n$,

\begin{equation}
\abs{2\sum_{j=1}^{2^n} \s_j^n (g) - \sum_{j=1}^{2^{n+1}}\T_j^n (g)}
\leq A(n+1)^{1/2}2^{n/2}.
\end{equation}

Let $G$ be the disjoint union of the sets $G_n$
and for each $n \ge 0$ and $1 \le j \le 2^n$ define 
$e_j^n : G \to \mathbb{C}$ via:

\begin{displaymath}
e_j^n(g) = \left\{ \begin{array}{ccc} 
\T_j^{n-1}(g), & \textrm{if $g \in G_{n-1}, n\ge 1$}\\
\varepsilon_j^n\s_j^n(g), & \textrm{if $g \in G_n$}\\
0, & \textrm{otherwise}
\end{array} \right.
\end{displaymath}
where   $\varepsilon_j^n = \pm 1$ is a   choice of signs for which the
inequality (2.5) below is satisfied.

To define $E$ let, as above,  $k_n = 3\cdot 2^n$ and let
$(p_n)_{n=0}^\infty$, $2 < p_n \leq 3$, be a decreasing sequence
converging  to $2$. The appropriate rate of decrease of the sequence
$(p_n)_{n=0}^\infty$ will be chosen later. 

Let $Z = \left(\sum_{n=0}^{\infty}\ell_{p_n}^{k_n}\right)_{2}$ 
which in our setting is:

\begin{displaymath}
Z = \bigl\{ f:G \to \mathbb{C} \,\big|\,
\sum_{n=0}^\infty \bigl(\sum_{g\in G_n}\abs{f(g)}^{p_n} \bigr)^{2/p_n}
< \infty\bigr\}.
\end{displaymath}

Define $E$  to be the closed linear span in $Z$ of 
$\{e_j^n |\, n\ge 0,\,1
\le j \le 2^n\}$. To show that $E$ fails the approximation property one
proceeds as follows:

For $n \ge 0$ and $1 \le j \le 2^n$ define $\alpha_j^n \in E^*$ by

\begin{equation}
\alpha_j^n (f) = 3^{-1}2^{-n} \sum_{g\in G_n}
\varepsilon_j^n\s_j^n(g^{-1})f(g). 
\end{equation}

When $n \ge 1$ the expression above equals 

\begin{equation}
\alpha_j^n (f) = 3^{-1}2^{1-n} \sum_{g\in G_{n-1}}
\T_j^{n-1}(g^{-1})f(g). 
\end{equation}

This follows from the fact that $\alpha_j^n(e_i^k)
=\delta_{ij}\cdot\delta_{kn}$ (because of the orthogonality of the
characters of a group) and then a linearity and continuity argument
shows that (2.2) and (2.3) agree on $E$.

Now let $B(E)$ be the space of bounded, linear operators on $E$, and
for each $n \ge 0$ define $\beta^n$ in the dual space $B(E)^*$ as:

\begin{displaymath}
\beta^n(T) = 2^{-n} \sum_{n=1}^{2^n}\alpha_j^n(T(e_j^n)),
\quad T\in B(E)
\end{displaymath}
Using $(2.2)$ we can rewrite $\beta^n$ as:

\begin{displaymath}
\beta^n(T) = 3^{-1}4^{-n}\sum_{g\in G_n} T
\biggl(\sum_{j=1}^{2^n}\varepsilon_j^n\s_j^n(g^{-1})e_j^n\biggr)(g)
\end{displaymath}
and from $(2.3)$ we get:

\begin{displaymath}
\beta^{n+1}(T) = 6^{-1}4^{-n}\sum_{g\in
G_n}T\biggl(\sum_{j=1}^{2^{n+1}}\T_j^n(g^{-1})e_j^{n+1}\biggr)(g)
\end{displaymath}
hence,

\begin{equation}
\beta^{n+1}(T) - \beta^n(T) =
3^{-1}2^{-n}\sum_{g\in G_n}T(\Phi_g^n)(g) 
\end{equation}
where,

\begin{displaymath}
\Phi_g^n = 2^{-n-1}\sum_{j=1}^{2^{n+1}}\T_j^n(g^{-1})e_j^{n+1} -
 2^{-n} \sum_{j=1}^{2^n} \varepsilon_j^n\s_j^n(g^{-1})e_j^n, \qquad g
\in G_n 
\end{displaymath}

Note that $\Phi_g^n \in E$ for every $g\in G_n$ and $n \geq 1$.
Now we estimate the right hand side of $(2.4)$. If $n \geq 1$ and $g \in
G_n$
then, 

\begin{displaymath}
 3^{-1}2^{-n} \sum_{g\in G_n}\abs{T(\Phi_g^n)(g)} 
\leq  \sup_{g \in G_n} \{\|T(\Phi_g^n)\|_{\infty}\} 
\leq \sup_{g \in G_n} \{\|T(\Phi_g^n)\|_{Z}\}. 
\end{displaymath}
Therefore,

\begin{displaymath}
\abs{\beta^{n+1}(T) - \beta^n(T)} \leq 
\sup_{g \in G_n} \{\|T(\Phi_g^n)\|_{Z}\} \text{ for every  $T \in
B(E)$.}
\end{displaymath} 

Note that from (2.1) we have that $|\Phi_g^n(h)| \leq
A(n+1)^{1/2}2^{-n/2}$  for $g,h \in G_n$. By applying lemma $(a)$ in
[D], the signs 
$\varepsilon_j^n$, $1 \le j \le 2^n$, can be chosen so that 

\begin{equation}
\abs{\Phi_g^n(h)} \leq A_2(n+1)^{1/2}2^{-n/2} \quad \text{for $g
\in G_n, h \in G_{n-1} \quad (n \geq 1)$}
\end{equation}
where $A_2$ is some absolute constant. An algebraic
argument shows that a similar estimate can be obtained for $g \in G_n$
and $h \in G_{n+1}$. In  brief, we have that there is an absolute
constant, say
$A$, such that,

\begin{equation}
\abs{\Phi_g^n(h)} \leq A(n+1)^{1/2}2^{-n/2} \quad \text{ for $g \in G_n$ 
 and  $h \in G_{n-1}\sqcup G_n\sqcup G_{n+1}$}.  
\end{equation}

Now, if $n \geq 1$ and $g \in G_n$ then,

\begin{eqnarray}
\| \Phi_g^n \|_{Z}^2 &=& \sum_{j=n-1}^{n+1}\biggl(\sum_{h\in
G_{j}}\abs{\Phi_g^n (h)}^{p_j}\biggr)^{2/p_j} \nonumber\\ &\leq & A^2
(n+1) 2^{-n}\bigl(
(3\cdot 2^{n-1})^{2/p_{n-1}} \,+\, (3\cdot 2^n)^{2/p_n} \,+\,
(3\cdot2^{n+1})^{2/p_{n+1}}\bigr) \nonumber\\ &\leq & 3A^2 (n+1)
2^{-n}\bigl(2^{2(n-1)/p_{n-1}} \,+\, 2^{2n/p_n} \,+\, 
2^{2(n+1)/p_{n+1}}
\bigr) \nonumber\\ &\leq & 18A^2 (n+1) 2^{2n(1/p_{n+1} - 1/2)} \nonumber
\end{eqnarray}
thus, 

\begin{equation}
\|\Phi_g^n \|_{Z} \leq 3\sqrt{2}A(n+1)^{1/2} 2^{n(1/p_{n+1} -
1/2)}
\end{equation}

Consider the set 
\begin{displaymath}
\mathcal{C} = \{e_1^0\}\cup\{(n+1)^2 \Phi_g^n | g\in G_n, n\geq 1\}
\end{displaymath}

The estimate in (2.7) clearly shows that when 
\begin{equation}
(n+1)^{5/2} 2^{n(1/p_{n+1} - 1/2)} \to 0 
\end{equation}
the set $\mathcal{C}$ becomes a relatively compact subset of $E$.
Obviously there are many choices for $(p_n)_{n=0}^{\infty}$,
$p_n\downarrow 2$, that satisfy $(2.8)$; in particular,  
  $1/p_n = 1/2 - 1/(n+1)^{\alpha}$ for any 
$\alpha < 1$ gives a sequence satisfying $(2.8)$. When $n(1/p_{n+1} -
1/2) = -3 \log_2 (n+1)$,  the sequence $(p_n)$ is the one with the
slowest (up
to a constant) possible rate of decrease for this construction. This
makes the space
$Z$ ``almost'' a weak Hilbert space in the sense that for some $K$,
$H_Z(m,K)\le m^{\log
\log m}$ for large $m$. Indeed, consider $F$, a subspace of 
$\left(\sum_{j=n+1}^{\infty}\ell_{p_j}^{k_j}\right)_{2}$ 
of dimension $m=m(k_1 + \cdots + k_n)$, where $m$ is the largest integer
 such that
$0 < 1/2 - 1/p_{n+1} < 1/\log_2(m)$. Then, 
$d(F,\ell_2^m) \leq T_2(F)C_2(F)$. The type 2 constant of $F$ is bounded
by an absolute constant independent of $m$, say $c_1$. The cotype 2
constant of $F$ can be estimated, using Tomczak's lemma, by the
cotype 2 constant
$C_{2,m}(\cdot)$ on
$m$ vectors  (see Section 5.25 in [T-J]):

\begin{eqnarray}
C_{2}(F)\le \sqrt{2} C_{2,m}(F) &\leq & \sqrt{2}  C_{p_{n+1}}(F)
m^{1/2 - 1/p_{n+1}} 
\nonumber\\
&\leq & \sqrt{2}  c_2 m^{1/2 - 1/p_{n+1}} 
\nonumber\\
&\leq & 2\sqrt{2}  c_2 \qquad\text{ (by the choice of  $m$)}\nonumber.
\end{eqnarray}

Hence, for $K:=2\sqrt{2}  c_1 c_2$ we obtain that $d(F,\ell_2^m) \leq
K$. 

Finally, to show that $E$ fails the approximation property the argument
in \cite{D} finishes as follows: for every $T \in B(E)$,

\begin{displaymath}
\abs{\beta^{n+1}(T) - \beta^n(T)} \leq 
\sup_{g \in G_n} \{\|T(\Phi_g^n)\|_{Z}\} \leq
(n+1)^{-2} \sup_{x\in \mathcal{C}} \|Tx\|_Z
\end{displaymath}
Also, 
\begin{displaymath}
\abs{\beta^0(T)} \leq \|Te_0^1\| \leq \sup_{x\in \mathcal{C}} \|Tx\|_Z
\end{displaymath}

Hence $\beta(T) = \lim_{n \to \infty} \beta^n(T)$ exists for all $T \in
B(E)$  and satisfies 
\begin{displaymath}
\abs{\beta(T)} \leq 3 \sup_{x\in \mathcal{C}} \|Tx\|_Z
\end{displaymath}

In particular, when $\mathcal{C}$ is compact, $\beta$ is a continuous
linear functional on $B(E)$ when $B(E)$ is given the topology of uniform
convergence on compact sets.

If $I_E$ is the identity map on $E$, it follows from the definition of
$\beta^n$ that $\beta^n(I_E) = 1$ for all
$n$, so $\beta(I_E) = 1$. On the other hand it is easy to see that
$\beta$
vanishes on the set of finite rank operators on $E$, thus $E$ cannot
have the approximation property.

\begin{remark} For $(p_n, k_n)_{n=0}^\infty$ as above, we
obtained an asymptotically Hilbertian space $Z$ which has a subspace
failing the approximation  property. The space $Z$ can be decomposed as
the direct sum of two subspaces,  say $Z_1$ and $Z_2$, all of whose
subspaces have the approximation property. Indeed, as in example 1.g.7
in \cite{LT}, it is enough to construct a subsequence $(p_{n_j})$ of
$(p_n)$ as follows: set $p_{n_1} = p_0$ and $k_{n_1} = k_0$. Having
chosen 
$p_{n_1} \cdots p_{n_j}$ (and their respective $k_{n_1} \cdots
k_{n_j}$),
choose $p_{n_{j+1}}$ such that if $F \subset \ell_p$  
$(2 < p < p_{n_{j+1}})$, has dimension 
$m \leq 2\cdot 5^{\sum_{i=1}^j k_{n_i}}$ then 
$m^{1/2 -1/p_{n_{j+1}}} \leq 2$  (in particular 
$d (F, \ell_2^m) \leq 2)$.
Now set $N_1= \{j  | p_{n_{2k+1}} \leq j < p_{n_{2(k+1)}}, k\geq 0\}$
and
$N_2 = \mathbb{N}-N_1$. Let 
$Z_i = \left(\sum_{j\in N_i}\ell_{p_j}^{k_j}\right)_{2}, i=1,2$. 
\end{remark}
\bs

Our example is best possible in light of current theory and the
current wisdom in the field.  First, it follows from the arguments in
\cite{J} that the spaces $Z_{1},Z_{2}$ have the property that every
subspace of every quotient space of $Z_{i}$, $i=1,2$ has a
decomposition of the form $Z_{i} = \left( \sum_{k=1}^{\infty} 
E_{k}\right )_{{\ell}_{2}}$,
where dim $E_{k}<\infty$  for each $k=1,2,3,\ldots $.  Also, the
argument
of Szarek \cite{S} shows that the spaces $Z_{i}$ have subspaces without
bases.
One might try to refine this example to produce $Z_{i}$'s for which
every subspace has a basis.  However, this may not be possible since it
is an open question whether Banach spaces for which every subspace has
a basis must be weak Hilbert.  Since the direct sum of weak Hilbert
spaces is weak Hilbert, and every subspace of a weak Hilbert space has
the approximation property, a positive answer to this question would
show that our construction cannot be improved to produce $Z_{i}$'s for
which every subspace has a basis.  
It was shown by Maurey and Pisier (see \cite{M}) that every separable
weak Hilbert space  $X$ has a finite dimensional decomposition.  That
is, there is a sequence of finite dimensional subspaces $E_{i}$ of $X$
so that for every $x\in X$ there is a unique sequence $x_{i}\in E_{i}$
so that $x = \sum_{i}x_{i}$.  However, it is an open 
question whether a separable weak
Hilbert space  must have a basis.  

\bibliographystyle{amsalpha}

\end{document}